\documentclass[oneside,english]{amsart}
\usepackage[T1]{fontenc}
\usepackage[latin1]{inputenc}
\usepackage{amssymb}
\usepackage{mathrsfs, tabularx}

\makeatletter
\theoremstyle{plain}
\newtheorem{thm}{Theorem}[section]
\newtheorem{prop}[thm]{Proposition}

\theoremstyle{definition}
\newtheorem{definition}[thm]{Definition}
\theoremstyle{remark}

\newtheorem*{remark*}{Remark}
\newtheorem*{acknowledgement*}{Acknowledgements}

\newcommand{\vfield}[1]{{\mathfrak X}\left( #1\right)}

\newcommand{\slaz}[0]{\setminus \{0\}}

\newcommand {\Set}[1] {\mathbb{#1}}
\newcommand{\setR}[0]{\Set{R}}
\newcommand{\setF}[0]{\Set{F}}
\newcommand{\setC}[0]{\Set{C}}

\def\cL{\mathcal{L}}

\usepackage{babel}
\makeatother
\begin{document}

\title[Semi-basic $1$-forms and  Helmholtz conditions]
{Semi-basic $1$-forms and  Helmholtz conditions for the inverse
problem of the calculus of variations}
\author[Bucataru]{Ioan Bucataru}
\address{Ioan Bucataru, Faculty of Mathematics, Al.I.Cuza
University, B-dul Carol 11, Iasi, 700506, Romania}
\urladdr{http://www.math.uaic.ro/\textasciitilde{}bucataru/}

\author[Dahl]{Matias F. Dahl}
\address{Matias F. Dahl, Institute of Mathematics, P.O.Box 1100, 02015
Helsinki University of Technology, Finland}
\urladdr{http://www.math.tkk.fi/\textasciitilde{}fdahl/}

\date{\today}

\begin{abstract}
We use Fr\"olicher-Nijenhuis theory to obtain global Helmholtz
conditions, expressed in terms of a semi-basic $1$-form, that
characterize when a semispray is locally Lagrangian. We also
discuss the relation between these Helmholtz conditions and their
classic formulation written using a multiplier matrix. When the
semi-basic $1$-form is $1$-homogeneous ($0$-homogeneous) we show
that two (one) of the Helmholtz conditions are consequences of the
other ones. These two special cases correspond to two inverse
problems in the calculus of variation: Finsler metrizability for a
spray, and projective metrizability for a spray.
\end{abstract}

\maketitle



\section{Introduction}
The inverse problem of the calculus of variations can be
formulated as follows. Under what conditions a system of second
order differential equations (SODE), on a $n$-dimensional manifold
$M$,
\begin{eqnarray} \frac{d^2x^i}{dt^2}+2G^i\left(x,\dot{x}\right)=0,
i\in \{1,2,...,n\}, \label{sode}
\end{eqnarray} can be derived from a variational principle? An
approach to this problem uses the Helmholtz conditions, which are
necessary and sufficient conditions for the existence of a
multiplier matrix $g_{ij}(x,\dot{x})$ such that
\begin{eqnarray}
g_{ij}(x,\dot{x})\left(\frac{d^2x^j}{dt^2}+2G^j\left(x,\dot{x}\right)\right)=
\frac{d}{dt}\left(\frac{\partial L}{\partial \dot{x}^i}\right) -
\frac{\partial L}{\partial {x}^i}, \label{variational_sode}
\end{eqnarray}
for some Lagrangian function $L(x,\dot{x})$. The multiplier matrix
$g_{ij}$ induces a symmetric $(0,2)$-type tensor field $g$ along
the tangent bundle projection. Geometric formulation of Helmholtz
conditions in terms of $g_{ij}$ were obtained by Sarlet
\cite{sarlet82} and expressed later using the tensor $g$ by
Martinez et al. \cite{martinez93a}. There are various approaches
to derive the Helmholtz conditions in both autonomous and
nonautonomous case. For discussions, see Crampin \cite{crampin81},
Krupkova and Prince \cite{krupkova07}, Morandi et al.
\cite{morandi90}.

In this paper we will study the inverse problem of calculus of
variations when the system of SODE in equation \eqref{sode} arise
from a semispray. In Theorem \ref{thm_Helm} we give a global
formulation for the Helmholtz conditions in terms of a semi-basic
$1$-form. If there exists a semi-basic $1$-form that satisfies
these Helmholtz conditions, the $1$-form is the Poincar\'e-Cartan
$1$-form of a locally defined Lagrangian function. Then the
original semispray is an Euler-Lagrange vector field for this
Lagrangian. In Section \ref{helm_mm}, we explain how these
Helmholtz conditions for a $1$-form correspond to the classic
formulation of Helmholtz conditions in terms of a multiplier
matrix. To derive the Helmholtz conditions in Theorem
\ref{thm_Helm} we use Fr\"olicher-Nijenhuis theory on $TM\slaz$
and geometric structures on $TM\slaz$ induced by the semispray.
See Sections 2 and 3, respectively.

It has been shown recently that for the case of Finsler spaces one
of the Helmholtz condition is a consequence of the other ones,
Prince \cite{prince08}. In \cite{sarlet07}, Sarlet claims that
this Helmholtz condition is redundant for homogeneity of any
order. In Theorem \ref{exactlst}, we prove that, depending on the
degree of homogeneity, one or two of the Helmholtz conditions can
be derived from the other ones. Therefore, in Section
\ref{lagrangianspray}, we show that a spray $S$ is Lagrangian if
and only if only if only two or three of the four Helmholtz
conditions are satisfied, depending on the degree of homogeneity.
In particular we discuss Helmholtz conditions for two important
inverse problems: projective metrizability and Finsler
metrizability.

For the projective metrizability of a spray $S$, we show that $S$
is an Euler-Lagrange vector field for a $1$-homogenous Lagrangian
if and only if two of the Helmholtz conditions, expressed in terms
of a semi-basic, $0$-homogeneous $1$-form, are satisfied. In
Section \ref{lagrangianspray}, we explain how these two Helmholtz
conditions correspond to the Rapcs\'ak conditions that
characterize projective metrizability \cite{rapcsak62}. For other
characterizations of projective metrizability of a spray, see
Klein \cite{klein62},  Klein and Voutier \cite{klein68}, Shen
\cite{shen01}, and Szilasi \cite{szilasi07}. For the case of a
flat spray the two Helmholtz conditions lead to Hamel's equations
studied recently by Crampin \cite{crampin08} and Szilasi
\cite{szilasi07}.

For $k>1$, we show that a spray $S$ is an Euler-Lagrange vector
field of a $k$-homogeneous Lagrangian if and only if three of the
Helmholtz conditions are satisfied. In particular, when $k=2$, we
obtain three Helmholtz conditions, expressed in terms of a
semi-basic, $1$-homogeneous $1$-form, for a spray $S$ to be
Finsler metrizable. In Section \ref{lagrangianspray}, we explain
how these three Helmholtz conditions are related to previous
discussions for the Finsler metrizability of a spray. See the work
of Crampin \cite{crampin07}, Krupka and Sattarov \cite{krupka85},
Muzsnay \cite{muzsnay06}, Prince \cite{prince08}, Szilasi and
Vattam\'ani \cite{szilasi02}.

An important tool in this work is the dynamical covariant
derivative induced by a semispray $S$. The notion of dynamical
covariant derivative was first introduced by Carin\~ena and
Martinez in \cite{carinena91} as a covariant derivative along the
tangent bundle projection. A recent discussion of various
connections associated to a semispray and their relation with the
dynamical covariant derivative is due to Sarlet \cite{sarlet07}.
See also \cite{crampin96, martinez93a}. Since all the geometric
structures that can be derived from a semispray $S$ are naturally
defined on the tangent bundle $TM$, we introduce, in Section
\ref{sec:DCD}, the dynamical covariant derivative as a tensor
derivation on $TM$ and study commutation formulae with geometric
structures induced by $S$. For a semispray, the dynamical
covariant derivative preserves the induced horizontal and vertical
distributions and hence will preserve semi-basic (vector valued)
forms. The restriction to semi-basic forms of the dynamical
covariant derivative coincides with the semi-basic derivation
studied by Grifone and Muzsnay \cite{grifone00}.

\section{Preliminaries} \label{preliminaries}
By a manifold $M$ we mean a second countable Hausdorff space that
is locally homeomorphic to $\setR^n$ with $C^\infty$-smooth
transition maps. Here $n\ge 1$ is the dimension of $M$. By $TM$ we
mean the tangent bundle $(TM, \pi, M)$ and by $TM\setminus\{0\}$
the tangent bundle with the zero section removed. The canonical
submersion $\pi:TM \to M$ induces a \emph{natural foliation} on
$TM$, whose leafs are tangent spaces $T_pM=\pi^{-1}(p)$, $p\in M$.
Local coordinates on $M$ will be denoted by $x^i$, while induced
coordinates on $TM$ will be denoted by $x^i, y^i$. Then $x^i$ are
transverse coordinates for the natural foliation, and $y^i$ are
coordinates for the leafs of this foliation.

Throughout the paper we assume that all objects are
$C^\infty$-smooth where defined. The ring of smooth functions on a
manifold $M$ is denoted by $C^\infty(M)$, the $C^{\infty}$ module
of $k$-forms is denoted by $\Lambda^k(M)$, and the $C^{\infty}$
module of vector fields is denoted by $\vfield{M}$. The
$C^{\infty}$ module of $(r,s)$-type tensor fields on $M$ is
denoted by $\mathcal{T}^r_s(M)$ and $\mathcal{T}(M)$ denotes the
tensor algebra on $M$.

By a \emph{vector valued $l$-form} ($l\ge 0$) on a manifold $M$ we
mean a $(1,l)$-type tensor field on $M$ that is anti-symmetric in
its $l$-arguments.

If $c\colon I\to M$, $c=(x^i)$ is a curve, we denote by $c'$ its
tangent $c'\colon I\to TM$, $c'(t) = (x^i, \dot x^i)$.  A curve
$c$ is \emph{regular} if $c'(t)\in TM\slaz$ for all $t\in I$.

\subsection{Fr\"olicher-Nijenhuis theory on $TM\slaz$}
In this section we give a quick review of the
Fr\"olicher-Nijenhuis theory.  For systematic treatments, see the
original paper of Fr\"olicher and Nijenhuis \cite{frolicher56} and
the book of Kolar et al. \cite{komisl93}.  In this paper we apply
this theory on $TM\slaz$, following Grifone and Muzsnay
\cite{grifone00}, Klein and Voutier \cite{klein68}, de Le\'on and
Rodrigues \cite{deleon89}, and Szilasi \cite{szilasi03}.

Suppose $A$ is a vector valued $l$-form on $TM\slaz$, and $\alpha$
is a $k$-form on $TM\slaz$ where $l\ge 0$ and $k\ge 1$. Then the
\emph{inner product} of $A$ and $\alpha$ is the $(k+l-1)$-form
$i_A\alpha$ defined as
\begin{eqnarray} \label{iaalpha} & & i_A\alpha(X_1,\cdots, X_{k+l-1})=  \\
\nonumber &  & \frac{1}{l!(k-1)!}\sum_{\sigma\in S_{k+l-1}}
\operatorname{sign}(\sigma)\ \alpha\left(A(X_{\sigma(1)},
\cdots,X_{\sigma(l)}), X_{\sigma(l+1)},\cdots, X_{\sigma(k+l-1)}
\right), \end{eqnarray} where $X_1, \ldots, X_{k+l-1} \in
\vfield{TM\slaz}$, and $S_{p}$ is the permutation group of
elements $1,\ldots, p$. When $l=0$, $A$ is a vector field on
$TM\slaz$ and $i_A\alpha$ is the usual inner product of $k$-form
$\alpha$ with respect to a vector field $A$. When $l=1$, $A$ is a
$(1,1)$-type tensor field and $i_A\alpha$ is the $k$-form
\begin{eqnarray}
i_A\alpha(X_1,\cdots, X_k) = \sum^k_{i=1}\alpha(X_1,..., AX_i,
..., X_k). \label{i1alpha}
\end{eqnarray}
We also define $i_A \alpha = 0$ when $\alpha \in
\Lambda^0(TM\slaz)=C^\infty(TM\slaz)$ and $A$ is any vector valued
$l$-form on $TM\slaz$.

One can define an \emph{exterior inner product}
$\overline{\wedge}$ on the graded algebra of vector valued
differential forms on $TM\slaz$ using a similar formula as
\eqref{i1alpha}, \cite{frolicher56}. In this work we will need
only the exterior inner product of a vector valued $k$-form $A$
with a $(1,1)$-type tensor $B$. In this case we define
$B\overline{\wedge} A$ as the vector valued $k$-form
\begin{eqnarray} B\overline{\wedge} A(X_1,\cdots, X_k) =
\sum^k_{i=1}B(X_1,..., AX_i, ..., X_k). \label{i1b}
\end{eqnarray}

Let $A$ be a vector valued $l$-form on $TM\slaz$, where $l\ge 0$.
Then the \emph{exterior derivative} with respect to $A$ is the map
$d_A\colon \Lambda^k(TM\slaz)\to \Lambda^{k+l}(TM\slaz)$ for $k\ge
0$,
\begin{eqnarray} d_A=i_A \circ d - (-1)^{l-1} d\circ i_A. \label{da}
\end{eqnarray}
A $k$-form $\omega$ on $TM\slaz$ is called $d_A$-\emph{closed} if
$d_A\omega=0$ and $d_A$-\emph{exact} if there exists $\theta \in
\Lambda^{k-l}(TM\slaz)$ such that $\omega=d_A\theta$.

When $A\in \vfield{TM\slaz}$ (that is, when $l=0$) and $k\ge 0$,
we obtain $d_A=\cL_A$, where $\cL_A$ is the usual Lie derivative
$\mathcal{L}_A\colon \Lambda^k(TM\slaz)\to \Lambda^k(TM\slaz)$. In
this case equation \eqref{da} is \emph{Cartan's formula}.

If $A=\operatorname{Id}$, then $l=1$ and $d_{\operatorname{Id}}
=d$ since $i_{\operatorname{Id}}\alpha = k \alpha$ for $\alpha\in
\Lambda^k(TM\slaz)$.

Suppose $A$ and $B$ are vector valued forms on $TM\slaz$ of
degrees $l\ge 0$ and $k\ge 0$, respectively. Then, the
\emph{Fr\"olicher-Nijenhuis bracket} of $A$ and $B$ is the unique
vector valued $(k+l)$-form $[A,B]$ on $TM\slaz$ such that
\cite{frolicher56},
\begin{eqnarray}
\label{dab} d_{[A,B]} = d_A \circ d_B - (-1)^{kl} d_B\circ d_A.
\end{eqnarray}
When $A$ and $B$ are vector fields (that is, when $k=l=0$), then
Fr\"olicher-Nijenhuis bracket $[A,B]$ coincides with the usual Lie
bracket $[A,B]=\cL_AB$.

When $A$ and $B$ are $(1,1)$-type tensor fields (that is, when
$k=l=1$), Fr\"olicher-Nijenhuis bracket $[A,B]$ is the vector
valued $2$-form \cite[p. 73]{komisl93}
\begin{eqnarray}
\label{abxy} [A,B](X,Y) & = & [AX, BY] + [BX, AY] + (AB + BA)[X,Y] \\
\nonumber & & - A[X,BY] - B[X, AY] -A[BX, Y] - B[AX,Y].
\end{eqnarray}
In particular,
\begin{eqnarray}
\label{aaxy} \frac{1}{2}[A,A](X,Y) = [AX, AY] + A^2[X,Y] - A[X,AY]
-A[AX, Y].
\end{eqnarray}
For a $(1,1)$-type tensor field $A$, the vector valued $2$-form
$N_A=(1/2)[A,A]$ is called the \emph{Nijenhuis tensor} of $A$.

For a vector field $X$ in $\vfield{TM\slaz}$ and a $(1,1)$-type
tensor field $A$ on $TM\slaz$ the Fr\"olicher-Nijenhuis bracket
$[X,A]=\mathcal{L}_XA$ is the $(1,1)$-type tensor field on
$TM\slaz$ given by
\begin{eqnarray}
\mathcal{L}_XA= \mathcal{L}_X \circ A - A \circ \mathcal{L}_X.
\label{xa}
\end{eqnarray}
Next commutation formulae on $\Lambda^k(TM\slaz)$, $k\ge 0$, will
be used throughout the paper, \cite{grifone00}:
\begin{eqnarray} \label{iadb} i_Ad_B-d_Bi_A=d_{B\circ A}-i_{[A,B]},\\
\label{lxia} \mathcal{L}_Xi_A-i_A \mathcal{L}_X=i_{[X,A]}, \\
i_Xd_A+ d_Ai_X=\mathcal{L}_{AX}- i_{[X,A]},
 \label{ixda}
 \end{eqnarray}
for $(1,1)$-type tensor fields $A,B$ and a vector field $X$ on
$TM\slaz$. We will refer to formula \eqref{ixda} as to the
generalized Cartan's formula, since by taking
$A=\operatorname{Id}$, it reduces to the usual Cartan formula.

\subsection{Homogeneous objects}
Suppose $k$ is an integer.  Then a function $f\in
C^{\infty}(TM\setminus\{0\})$ is said to be \emph{positively
$k$-homogeneous} (or briefly \emph{$k$-homogeneous}) if $f(\lambda
y)=\lambda^kf(y)$ for all $\lambda>0$ and $y\in TM\slaz$.  By
Euler's theorem, a function $f\in C^{\infty}(TM\setminus\{0\})$ is
$k$-homogeneous if and only if $\mathcal{L}_{\mathbb{C}}f=kf$,
where $\setC\in \vfield{TM}$ is the \emph{Liouville vector field}
(or \emph{dilatation vector field}) defined as $\setC(y) =
(y+sy)'(0)$. In local coordinates $(x^i, y^i)$ for $TM$,
\begin{equation}
\mathbb{C}=y^i\frac{\partial}{\partial y^i}. \label{Liouville}
\end{equation}
Using vector field $\setC$, we also define homogeneity for other
objects on $TM\slaz$. A vector field $X\in
\vfield{TM\setminus\{0\}}$ is $k$-homogeneous if and only if
$\mathcal{L}_{\mathbb{C}}X=(k-1)X$. Alternatively, a vector field
is $k$-homogeneous if its flow is $k$-homogeneous. For example,
Liouville vector field $\setC$ is $1$-homogeneous.  A $p$-form
$\omega \in \Lambda^p(TM\setminus\{0\})$ is $k$-homogeneous if and
only if $\mathcal{L}_{\mathbb{C}}\omega=k\omega$. Lastly, a
$(1,1)$-tensor $L$ on $TM\slaz$ is $k$-homogeneous if and only if
$\mathcal{L}_{\mathbb{C}}L=(k-1)L$.

\subsection{Vertical calculus on $TM\slaz$}
Next, we define the canonical tangent structure $J$ on $TM\slaz$,
which is a $(1,1)$-type tensor on $TM\slaz$. Then, the
Fr\"olicher--Nijenhuis theory gives a particular differential
calculus with operators $i_J$ and $d_J$. These operators are well
suited for studying Finsler and Lagrange geometries on $TM\slaz$
\cite{BM07, grifone00,deleon89, morandi90,szilasi03}. They also
play a key role in this paper.

The \emph{vertical subbundle} is defined as
\begin{eqnarray}
  VTM = \{ \xi\in TTM : (D\pi)(\xi)=0\}.
\end{eqnarray}
Then the map $V_u: u\mapsto V_u=VTM\cap T_uTM$ defines the
\emph{vertical distribution} $V$. It is a $n$-dimensional,
integrable distribution, being tangent to the natural foliation.
That is, any vertical vector $u\in VTM$ can be written as $u =
(y+tz)'(0)$ for some vectors $y,z\in TM$ with $\pi(y)=\pi(z)$. An
important vertical vector field on $TM\slaz$ is the Liouville
vector field \eqref{Liouville}.

On $TM\slaz$ \emph{the tangent structure} (or the \emph{vertical
endomorphism}) is the $(1,1)$-type tensor $J$ defined as
$$J(\xi)=\left(\tau(\xi)+t(D\pi)(\xi)\right)'(0), \forall \xi \in
TTM. $$ Here $\tau:TTM \to TM$ is the canonical submersion of the
second order iterated tangent bundle. Locally,
\begin{eqnarray}
  J=\frac{\partial}{\partial y^i}\otimes dx^i. \label{tangent}
\end{eqnarray}
Tensor $J$ satisfies $J^2=0$ and
$\operatorname{Ker}J=\operatorname{Im}J=VTM$ and $J$ is
$0$-homogeneous since $\mathcal{L}_{\mathbb{C}}J = [{\mathbb{C}},
J]=-J$, \cite{grifone72}. An important notion in this work is that
of semi-basic forms, \cite{grifone00,deleon89}.

\begin{definition} \label{def:semi-basic} Consider $k\ge 1$.
\begin{itemize} \item[i)]
A $k$-form $\omega$ on $TM\slaz$ is called \emph{semi-basic} if
$\omega(X_1, ..., X_k)=0$, when one of the vectors $X_i$,
$i\in\{1,...,k\}$ is vertical.
\item[ii)] A vector valued $k$-form $A$ on $TM\slaz$ is called \emph{semi-basic}
if it takes values in the vertical bundle and $A(X_1, ...,
X_k)=0$, when one of the vectors $X_i$, $i\in\{1,...,k\}$ is
vertical.
\end{itemize} \end{definition}

If a $k$-form $\omega$ is semi-basic then using formula
\eqref{i1alpha} we obtain that $i_J\omega=0$. The converse is true
only if $k=1$. In other words, a $1$-form $\theta\in
\Lambda^1(TM\slaz)$ is semi-basic if and only if $i_J\theta=0$.
Moreover, any semi-basic $1$-form $\theta\in \Lambda^1(TM\slaz)$
can be written as $\theta=i_J\omega$, for a (non unique) $1$-form
$\omega \in \Lambda^1(TM\slaz)$. Semi-basic $1$-forms are
annihilators for the vertical distribution. In local coordinates
$(x^i,y^i)$ for $TM\slaz$, a semi-basic $1$-form can be expressed
as $\theta=\theta_i(x,y)dx^i$.

If a vector valued $k$-form $A$ is semi-basic then $J\circ A=0$
and $A\overline{\wedge}J=0$. The converse is true only if $k=1$. A
vector valued $1$-form $A$ on $TM\slaz$ is semi-basic if and only
if $J\circ A=0$ and $A\circ J=0$. It follows that the tangent
structure $J$ is a vector valued, semi-basic $1$-form, and its
Nijenhuis tensor $N_J=(1/2)[J,J]$ vanishes. Hence equation
\eqref{dab} implies that
\begin{eqnarray}
d_J^2=d_J\circ d_J=0. \label{dj2} \end{eqnarray} Formula
\eqref{dj2} shows that any $d_J$-exact form on $TM\slaz$ is also
$d_J$-closed.

For semi-basic forms, $d_J$ is the exterior differential along the
leafs of the natural foliation and from formula \eqref{dj2} it
satisfies a local Poincar\'e lemma, \cite{vaisman73, vaisman88}.
Therefore, $d_J$-closed semi-basic forms on $TM\slaz$ are locally
$d_J$-exact. Note that a local Poincar\'e lemma does not hold true
if we do not restrict $d_J$ to semi-basic forms, as it has been
pointed out in \cite[p.173]{morandi90}. Locally, a semi-basic
$1$-form $\theta =\theta_i dx^i$ is $d_J$-closed and hence locally
$d_J$-exact if and only if the matrix $(\partial \theta_i/\partial
y^j)$ is symmetric. The relation between $d_J$-closed and
$d_J$-exact semi-basic forms on $TM\setminus\{0\}$ has been
discussed by Klein \cite{klein62} for the homogeneous case. In
this context it has been shown by Klein and Voutier \cite{klein68}
and de Le\'on and Rodrigues \cite{deleon89} that a semi-basic
$p$-form, $k$-homogeneous with $p\neq -k$, is $d_J$-closed if and
only if it is $d_J$-exact. In Proposition \ref{lem_hil} we will
specialize this result to homogeneous, semi-basic $1$-forms.

\begin{definition} A semi-basic $1$-form $\theta \in
\Lambda^1(TM\slaz)$ is called \emph{non-degenerate} if $d\theta$
is a symplectic form on $TM\slaz$. \end{definition} A semi-basic
$1$-form $\theta =\theta_i dx^i$ is non-degenerate if and only if
the matrix with entries $(\partial \theta_i/\partial y^j)$ is
non-degenerate.

\section{Semisprays and nonlinear connections} \label{semispray}

A system of second order differential equations (SODE) on a
manifold $M$, whose coefficients functions do not depend
explicitly on time, can be viewed as a special vector field on
$TM\slaz$, which is called a semispray. If the coefficients
functions of the SODE are $2$-homogeneous functions, then the
corresponding vector field is called a spray. In the affine
context, the notion of spray was introduced by Ambrose et al.
\cite{ambrose60} and later extended by Dazord \cite{dazord66}.

In this section, we start with a semispray $S$ and consider
induced geometric structures that will be useful to express
necessary and sufficient conditions for $S$ to be Lagrangian.
These geometric structures are defined using the nonlinear
connection induced by a semispray, which was considered first by
Crampin \cite{crampin71} and Grifone \cite{grifone72}. A nonlinear
connection can be characterized using horizontal and vertical
projectors, horizontal lifts, almost product structures or almost
complex structures, see \cite{BM07,grifone00,deleon89,miron94,
morandi90, szilasi03}. We point out some important features of the
induced geometric objects in the homogeneous case that will be
used in the paper.

\subsection{Semisprays and nonlinear connections} \label{snc}
\begin{definition} \label{def_spray} \begin{itemize}
\item[i)] A \emph{semispray} (or a \emph{second order vector field}) on
$M$ is a vector field $S\in \vfield{TM\slaz}$ such that
$JS=\mathbb{C}$. \item[ii)] A \emph{spray} on $M$ is a semispray
$S$ that is $2$-homogeneous as a vector field.
\end{itemize}
\end{definition} Locally, a semispray $S$ on $M$ can be written as
\begin{eqnarray}
\label{spray} S=y^i\frac{\partial}{\partial x^i} - 2
G^i(x,y)\frac{\partial}{\partial y^i},
\end{eqnarray}
for some functions $G^i$ called \emph{semispray coefficients} of
$S$. Functions $G^i$ are defined on domains of induced coordinate
charts on $TM\slaz$.

A spray on $M$ is a vector field $S\in \vfield{TM\slaz}$ such that
$JS=\mathbb{C}$ and $[\mathbb{C}, S]=S$. For a spray $S$ functions
$G^i$ in formula \eqref{spray} are $2$-homogeneous functions where
defined.
\begin{definition} \label{def_geod}
A regular curve $c: I \to M$ is a \emph{geodesic} of a semispray
$S$ if $S\circ c'=c''$.
\end{definition}
If $c(t)=(x^i(t))$ is a regular curve on $M$, then $c$ is a
geodesic of semispray $S$ in equation \eqref{spray} if it
satisfies the system of second order ordinary differential
equations
\begin{equation}
\label{eq_geod}
  \frac{d^2x^i}{dt^2}+ 2G^i\left(x, \frac{dx}{dt}\right)=0.
\end{equation}
Next we consider some tensors on $TM\slaz$ associated with a
semispray: horizontal and vertical projections $h$ and $v$, almost
product and complex structures $\Gamma$ and $\setF$, the Jacobi
endomorphism $\Phi$ and the curvature tensor $R$.

\begin{definition} \label{def_non} A \emph{nonlinear connection} (or a \emph{horizontal distribution}) on $M$
is defined by an $n$-dimensional distribution $H: u\in TM\slaz \to
H_u\subset T_u(TM\slaz)$ that is supplementary to the vertical
distribution $V$, which means that $T_u(TM\slaz)=H_u\oplus V_u$,
for all $u\in TM\slaz$.
\end{definition}

The \emph{horizontal projector} $h$ and \emph{vertical projector}
$v$ are $(1,1)$-type tensors on $TM\slaz$ defined as
\cite{grifone72},
\begin{equation}
h=\frac{1}{2}\left(\operatorname{Id}-\mathcal{L}_SJ\right), \quad
v=\frac{1}{2}\left(\operatorname{Id}+\mathcal{L}_SJ\right).
\label{hproj}
\end{equation}
Locally,
\begin{eqnarray*}
\label{hvloc} h = \frac{\delta}{\delta x^i}\otimes dx^i, \quad v =
\frac{\partial}{\partial y^i}\otimes \delta y^i,
\end{eqnarray*}
where
\begin{equation}
\frac{\delta}{\delta x^i}=\frac{\partial}{\partial x^i}-N^j_i
\frac{\partial}{\partial y^j}, \quad \delta y^i=dy^i+ N^i_jdx^j, \
\mbox{ and }\ N^i_j=\frac{\partial G^i}{\partial y^j}.
\label{delta}
\end{equation}
Functions $N^i_j$ are called the \emph{nonlinear coefficients}
associated to semispray $S$. The $(1,1)$-type tensor field
\begin{equation}
\Gamma=-\mathcal{L}_SJ \label{lsj} \end{equation} used to define
the horizontal and vertical projectors in formulae \eqref{hproj}
is called the \emph{almost product structure} induced by semispray
$S$, \cite{grifone72}. It can be written as $\Gamma=h-v$ and
therefore $\Gamma^2=\operatorname{Id}$.

The \emph{almost complex structure} is the (1,1)-type tensor field
on $TM\slaz$ given by \cite{grifone00,morandi90}
\begin{equation}
\mathbb{F}=h\circ\mathcal{L}_Sh-J. \label{complex} \end{equation}
Locally,
\begin{equation}
\mathbb{F}=\frac{\delta}{\delta x^i}\otimes \delta y^i -
\frac{\partial}{\partial y^i}\otimes dx^i. \label{def_complex}
\end{equation} It follows immediately that $\mathbb{F}^2 = -\operatorname{Id}$.
Moreover, the following formulae for the above considered
$(1,1)$-type tensor fields will be useful throughout the paper:
$$\mathbb{F}\circ J=h, \quad J\circ \mathbb{F}=v, \quad v\circ
\mathbb{F}=\mathbb{F}\circ h=-J, \quad h\circ \mathbb{F}
=\mathbb{F}\circ v = \mathbb{F}+ J.$$

The \emph{Jacobi endomorphism} $\Phi$ is defined as the
$(1,1)$-type tensor field
\begin{eqnarray}
\Phi=v\circ \mathcal{L}_Sh = -v\circ \mathcal{L}_Sv.
\label{jacobiphi}
\end{eqnarray} The Jacobi endomorphism $\Phi$ is a semi-basic vector valued
$1$-form and it is also called the \emph{Douglas tensor}
\cite{grifone00}. Jacobi endomorphism $\Phi$ has been defined as a
$(1,1)$-type tensor field along the tangent bundle projection in
\cite{carinena91,crampin96,martinez93a}. Locally,
\begin{eqnarray}
\Phi= R^i_j \frac{\partial}{\partial y^i} \otimes dx^j,
\end{eqnarray}
where
\begin{equation}
R^i_j = 2\frac{\partial G^i}{\partial x^j} - S\left(\frac{\partial
G^i}{\partial y^j}\right) - \frac{\partial G^i}{\partial
y^r}\frac{\partial G^r}{\partial y^j}. \label{jacobi_end}
\end{equation}
The Jacobi endomorphism $\Phi$ has been used to study various
aspects of an SODE: variational equations \cite{BM07,carinena91},
symmetries \cite{BM07,carinena91, martinez93}, separability
\cite{martinez93}, linearizability \cite{crampin96} as well as to
express one of the Helmholtz condition of the inverse problem of
the calculus of variation \cite{crampin81,
krupkova97,martinez93a,sarlet82, sarlet07}.

The \emph{curvature tensor} $R$ of a nonlinear connection $N$ is
defined as the Nijenhuis tensor of the horizontal projector $h$,
$R=(1/2)[h,h]$. Locally,
\begin{equation}
  R=R^k_{ij}dx^i\wedge dx^j \otimes \frac{\partial}{\partial y^k},
\label{curvature} \end{equation} where
\begin{eqnarray}
R^k_{ij} = \frac{\delta N^k_i}{\delta x^j} - \frac{\delta
N^k_j}{\delta x^i}. \label{rkij}
\end{eqnarray} For the curvature tensor $R$ we have that
$R(X,Y)=R(hX,hY)=v[hX,hY]$ for all $X, Y\in \vfield{TM\slaz}$.
Therefore $R$ is a semi-basic, vector valued $2$-form that
vanishes if and only if the horizontal distribution is integrable.
If the horizontal distribution is integrable, then it is tangent
to a foliation that is transverse to the natural foliation and
$d_h$ is the exterior differentiation along the leafs of this
transverse foliation. It follows that for an integrable horizontal
distribution we have that $d_h^2=d_R=0$ and the restriction of the
differential operator $d_h$ to forms tangent to the transverse
foliation satisfies a local Poincar\'e lemma, \cite{vaisman73}.
Consequently, for a flat nonlinear connection, $d_h$-exact
$1$-forms tangent to the transverse foliation are locally
$d_h$-closed.

The curvature tensor $R$ can be obtained directly from the Jacobi
endomorphism $\Phi$ through the following formula,
\cite{grifone00,martinez93, szilasi03}
\begin{eqnarray} 3[J,\Phi]+R=0. \label{jphir} \end{eqnarray}
One can also recover the Jacobi endomorphism $\Phi$ from the
curvature tensor $R$ through the following formula
\begin{eqnarray} \Phi=i_SR + v\circ \mathcal{L}_{vS} h. \label{phir}
\end{eqnarray}
Indeed for a vector field $X$ on $TM\slaz$, we have
$\Phi(X)=v[S,hX]$ and $R(S,X)=v[hS,hX]$. Therefore,
$\Phi(X)=R(S,X)+v[vS,hX]$, which proves formula \eqref{phir}.

If $S$ is a spray then by Euler's theorem, the nonlinear
coefficients $N^i_j$ defined by formula \eqref{delta} are
$1$-homogeneous. Using the homogeneity of a spray $S$ and the
horizontal projector \eqref{hproj} it follows that $S=hS$, which
implies that $S$ has the local expression
\begin{eqnarray}
\label{spraydelta}
  S = y^i \frac{\delta}{\delta x^i}.
\end{eqnarray}
Therefore, for a spray $S$, we have that $vS=0$ and formula
\eqref{phir} gives
\begin{eqnarray} \Phi=i_SR. \label{phir2}
\end{eqnarray} In local coordinates formula \eqref{phir2} can be
written as
\begin{eqnarray} R^i_j(x,y)=R^i_{kj}(x,y)y^k, \label{rry}
\end{eqnarray} and connects the Jacobi endomorphism $R^i_j$ given
by formula \eqref{jacobi_end} and the curvature tensor $R^i_{kj}$
given by formula \eqref{rkij}.

\subsection{Dynamical covariant derivative}
\label{sec:DCD}

When a semispray $S$ is given on a manifold $M$, the Lie
derivative $\mathcal{L}_S$ defines a tensor derivation on
$TM\slaz$. However, the derivation $\mathcal{L}_S$ does not
preserve the geometric structures introduced in Section \ref{snc}.
In this section we show how to modify the derivation
$\mathcal{L}_S$ to obtain a tensor derivation on $TM\slaz$ that
preserves these geometric structures. This derivation is called
the \emph{dynamical covariant derivative} of the semispray. The
notion of dynamical covariant derivative induced by a semispray
was first introduced by Carin\~ena and Martinez in
\cite{carinena91} as a derivation of degree $0$ along the tangent
bundle projection, see also \cite{crampin96, martinez93,
martinez93a,szilasi03}. It was also studied as a \emph{semi-basic
derivation} of semi-basic forms by Grifone and Muzsnay
\cite{grifone00}. An extensive discussion about the dynamical
covariant derivative $\nabla$ and other linear connection along
the tangent bundle projection, which are associated to a
semispray, is due to Sarlet \cite{sarlet07}.

\begin{definition} \label{deriv} A map $\nabla:
\mathcal{T}(TM\slaz) \to \mathcal{T}(TM\slaz)$ is said to be a
\emph{tensor derivation} on $TM\slaz$ if it satisfies the
following conditions:
\begin{itemize} \item[i)] $\nabla$ is $\mathbb{R}$-linear;
\item[ii)] $\nabla$ is type preserving, which means that
$\nabla(\mathcal{T}^r_s(TM\slaz))\subset
\mathcal{T}^r_s(TM\slaz)$, for each pair $(r,s)$ in
$\mathbb{N}\times \mathbb{N}$; \item[iii)] $\nabla$ obeys the
Leibnitz rule, which means that $\nabla(T\otimes S)=\nabla T
\otimes S + T \otimes \nabla S$ for any tensor fields $T, S$ on
$TM\slaz$;
\item[iv)] $\nabla$ commutes with any contractions.
\end{itemize} \end{definition}

For a semispray $S$ on $M$, let us consider the
$\mathbb{R}$-linear map $\nabla_0: \vfield{TM\slaz} \to
\vfield{TM\slaz}$
\begin{equation}
\nabla_0 X=h[S,hX] + v[S,vX], \forall X\in \vfield{TM\slaz}.
\label{def_nabla}
\end{equation}
One can immediately check that
\begin{equation} \nabla_0(fX)=S(f)\nabla_0 X + f \nabla_0 X, \forall f
\in C^{\infty}(TM\slaz), \forall X\in \vfield{TM\slaz}.
\label{nabla1}
\end{equation}
Any tensor derivation on $TM\slaz$ is completely determined by its
action on smooth functions and vector fields on $TM\slaz$, see
\cite[p. 1217]{szilasi03}. Therefore there exists a unique tensor
derivation $\nabla$ on $TM\slaz$ such that
$$ \nabla|_{C^{\infty}(TM\slaz)}=S, \quad
\nabla|_{\vfield{TM\slaz}}=\nabla_0.$$ We will call the tensor
derivation $\nabla$, the \emph{dynamical covariant derivative}
induced by the semispray $S$.

Next, we will obtain some alternative expressions for the action
of the dynamical covariant derivative $\nabla$ on vector fields,
forms and vector-valued forms on $TM\slaz$.

From formula \eqref{def_nabla}, we have that the action of
$\nabla$ on $\vfield{TM\slaz}$ can be written as
\begin{equation} \nabla=h\circ \mathcal{L}_S \circ h + v\circ \mathcal{L}_S \circ
v. \label{nabla2} \end{equation} Since $\mathcal{L}_S
h=\mathcal{L}_S\circ h - h\circ \mathcal{L}_S$, it follows that
formula \eqref{nabla2} can be written as
\begin{equation} \nabla=\mathcal{L}_S+ h\circ \mathcal{L}_S h + v\circ \mathcal{L}_Sv.
\label{nabla3} \end{equation} Formula \eqref{nabla3} can be
further expressed as
\begin{equation} \nabla=\mathcal{L}_S +\Psi, \label{nabla4} \end{equation}
where
\begin{equation} \Psi= h\circ \mathcal{L}_S h + v\circ
\mathcal{L}_S v= \Gamma \circ \mathcal{L}_S h = (\mathbb{F}+J) -
\Phi \label{def_psi}
\end{equation} is a (1,1)-type tensor field on $TM\slaz$.
Decomposition \eqref{nabla4} of the dynamical covariant derivative
$\nabla$ can be compared with decomposition formula (96) in
\cite{martinez93a}.

Let $\omega$ be a $k$-form on $TM\slaz$. Since $\nabla$ satisfies
the Leibnitz rule, we obtain
\begin{equation}
(\nabla \omega)(X_1,...,X_k) = \nabla(\omega(X_1,...,X_k)) -
\sum^k_{i=1}\omega(X_1,..., \nabla X_i, ..., X_l).
\label{nabla_omega}\end{equation} Using expressions
\eqref{nabla_omega} and \eqref{nabla4} we obtain that the
dynamical covariant derivative $\nabla$ has the following
expression on $\Lambda^k(TM\slaz)$
\begin{equation} \nabla = \mathcal{L}_S - i_{\Psi}.
\label{nabla_k} \end{equation}

The action of $\nabla$ on vector valued $k$-forms on $TM\slaz$ can
be defined using a formula similar with \eqref{nabla_omega}. We
obtain that for a vector valued $k$-form $A$ on $TM\slaz$, its
dynamical covariant derivative is given by
\begin{eqnarray}
\nabla A=\mathcal{L}_SA+ \Psi \circ A - A\overline{\wedge} \Psi.
\label{nablaa}
\end{eqnarray} Formula \eqref{nablaa} coincides with the semi-basic
derivation acting on semi-basic vector valued forms considered by
Grifone and Muzsnay \cite[Proposition 4.4]{grifone00}. When $k=1$
and $A$ is a $(1,1)$-type tensor field on $TM\slaz$, we obtain
that its dynamical covariant derivative is given by
\begin{equation} \nabla A= \mathcal{L}_SA + \Psi\circ A - A\circ \Psi.
\label{nabla_A} \end{equation}

Next theorem shows that the dynamical covariant derivative
$\nabla$ preserves by parallelism the geometric structures induced
by a semispray $S$.
\begin{thm} \label{prop:nabla} Consider $\nabla$ the dynamical
covariant derivative induced by a semispray $S$ and $k\geq 0$.
\begin{itemize}
\item[i)] $\nabla h=0$, $\nabla v=0$, which means that $\nabla$
preserves the horizontal and vertical distributions;
\item[ii)] $\nabla J=0$, $\nabla \mathbb{F}=0$, which means that
$\nabla$ acts identically on both vertical and horizontal
distributions (see also formulae \eqref{nabladelta} and
\eqref{nablapartial} below);
\item[iii)] The restriction of $\nabla$ to $\Lambda^k(TM\slaz)$ and the
exterior differential operator $d$ satisfies the commutation
formula
\begin{eqnarray}
d\nabla - \nabla d & =& d_{\Psi}. \label{dnabla}
\end{eqnarray}
\item[iv)] The restriction of $\nabla$ to $\Lambda^k(TM\slaz)$ satisfies
the following commutation rule: \begin{eqnarray} \label{ihnabla}
\nabla i_A-i_A \nabla & = & i_{\nabla A},
\end{eqnarray} for any  vector valued $l$-form $A$ on $TM\slaz$.
If $l=1$ and $A\in \{h,v, J, \Gamma, \mathbb{F}\}$ then
\begin{eqnarray} \nabla i_A-i_A \nabla & = & 0. \label{danabla0}
\end{eqnarray}
\end{itemize} \end{thm}
\begin{proof}
From formula \eqref{def_psi}, which defines the $(1,1)$-type
tensor $\Psi$, it follows that \begin{eqnarray} \label{hpsi}
h\circ \Psi -\Psi\circ h = \mathcal{L}_Sh, \\ \label{jpsi} J\circ
\Psi - \Psi\circ J = \mathcal{L}_SJ, \\
\label{fpsi} \mathbb{F}\circ \Psi - \Psi\circ \mathbb{F} =
\mathcal{L}_S\mathbb{F}.  \end{eqnarray} Using formula
\eqref{nabla_A}, we obtain that the first two items of the
proposition are true.

From formula \eqref{nabla_k} it follows that
$$d\nabla = d\mathcal{L}_S - di_{\Psi} = \mathcal{L}_S d -
i_{\Psi} d + d_{\Psi}=\nabla d + d_{\Psi}$$ and hence formula
\eqref{dnabla} is true.

We will mainly need formula \eqref{ihnabla} for $l=0$ or $l=1$. We
will prove it for $l=1$. Using formulae \eqref{nabla_k},
\eqref{lxia} and \eqref{nabla_A}, we have
$$ \nabla i_A-i_A \nabla=\mathcal{L}_S i_A- i_A\mathcal{L}_S -
i_{\Psi}i_A + i_Ai_{\Psi} = i_{[S,A]}-i_{A\circ \Psi}+ i_{\Psi
\circ A}=i_{\nabla A}. $$ Using first two items of the theorem and
formula \eqref{ihnabla} we obtain commutation formulae
\eqref{danabla0}.
\end{proof}
From Theorem \ref{prop:nabla} we obtain that $\nabla J=0$ and
$\nabla i_J=i_J\nabla$ and hence the dynamical covariant
derivative $\nabla$ preserves semi-basic (vector valued) forms.
The restriction of $\nabla$ to semi-basic forms coincides with the
\emph{semi-basic derivation} studied by Grifone and Muzsnay
\cite{grifone00}. Commutation rule \eqref{ihnabla} shows that the
dynamical covariant derivative $\nabla$ is a self-dual derivation
in the sense of \cite[Theorem 3.2]{martinez93a}.

To express the action of $\nabla$, let us first note that
$$ \left[S, \frac{\partial}{\partial y^i}\right] = -
\frac{\delta}{\delta x^i} + N^k_i \frac{\partial}{\partial y^k}, \
\left[S, \frac{\delta}{\delta x^i}\right] =
N^k_i\frac{\delta}{\delta x^k} + R^k_i \frac{\partial}{\partial
y^k}.$$ Therefore, it follows that
\begin{eqnarray} \label{nabladelta} \nabla \frac{\delta}{\delta
x^i}  = h\left[S, \frac{\delta}{\delta
x^i}\right] =  N^k_i \frac{\delta}{\delta x^k}, \\
\label{nablapartial} \nabla \frac{\partial}{\partial y^i} =
v\left[S, \frac{\partial}{\partial y^i}\right] =  N^k_i
\frac{\partial}{\partial y^k},
\end{eqnarray}
and hence $\nabla$ coincides with the covariant derivative studied
in \cite{bucataru07,BM07}. Since horizontal and vertical vector
fields can be projected onto vector fields along the tangent
bundle projection, one can also project formulae
\eqref{nabladelta} or \eqref{nablapartial} and obtain the
dynamical covariant derivative along the tangent bundle projection
studied in \cite{crampin96,krupkova07, martinez93,
martinez93a,szilasi03}.

The next proposition shows that when $S$ is a spray the dynamical
covariant derivative has more properties.
\begin{prop} \label{prop:nabla2} Consider $\nabla$ the dynamical
covariant derivative induced by a spray $S$.\begin{itemize}
\item[i)] $\nabla S=0$ and $\nabla \mathbb{C}=0$, \item[ii)] $\nabla
i_S=i_S\nabla$ and $\nabla
i_{\mathbb{C}}=i_{\mathbb{C}}\nabla$.\end{itemize} \end{prop}
\begin{proof}
Since $\Psi(S)=0$ and $\Psi(\mathbb{C})=S$ we obtain using formula
\eqref{nabla4} that $\nabla S=0$ and $\nabla\mathbb{C}=0$. Second
part follows from formula \eqref{ihnabla} for $l=0$ and $A\in \{S,
\mathbb{C}\}$.
\end{proof}

\section{Semi-basic $1$-forms and Helmholtz conditions}

In Section 5 we show that the geodesics of a semispray $S$ are
solutions of the Euler-Lagrange equations for some Lagrangian $L$
if and only if there exists a semi-basic $1$-form $\theta \in
\Lambda^1(TM\slaz)$ such that the $1$-form $\mathcal{L}_S\theta$
is closed. We first find necessary and sufficient conditions,
called Helmholtz conditions, for a semi-basic $1$-form $\theta \in
\Lambda^1(TM\slaz)$ such that the $1$-form $\mathcal{L}_S\theta$
is closed. We then relate these Helmholtz conditions with their
classic formulation in terms of a multiplier matrix. Finally, we
show that for a spray and a homogeneous, semi-basic $1$-form
$\theta \in \Lambda^1(TM\slaz)$, the $1$-form
$\mathcal{L}_S\theta$ is closed if and only if it is exact.
Moreover, depending on the degree of homogeneity, some of the
Helmholtz conditions can be derived from the other ones.

\subsection{Helmholtz conditions for semi-basic $1$-forms}

Next theorem provides necessary and sufficient conditions for a
semi-basic $1$-form $\theta\in \Lambda^1(TM\slaz)$ such that the
$1$-form $\mathcal{L}_S\theta$ is closed.

\begin{thm}
\label{thm_Helm} Let $S$ be a semispray on $M$ and let $\theta$ be
a semi-basic $1$-form on $TM\slaz$. Then $\mathcal{L}_S\theta$ is
closed if and only if it satisfies the following \emph{Helmholtz
conditions}
\begin{eqnarray}
  d_h\theta=0, \quad
  d_J\theta=0, \quad
  \nabla d\theta=0, \quad
  d_{\Phi}\theta=0. \label{Helmholtz}
\end{eqnarray}
\end{thm}

\begin{proof}
From formulae \eqref{nabla_k} and \eqref{def_psi} it follows that
for the $2$-form $d\theta$ we have
\begin{eqnarray}
\mathcal{L}_S d\theta= \nabla d\theta + i_{\mathbb{F}+J}d\theta -
d_{\Phi}\theta. \label{nabladt1}\end{eqnarray}

For a semi-basic $1$-form $\theta \in \Lambda^1(TM\slaz)$ we have
\begin{equation} (d\theta)(JX, JY)=(JX)((\theta\circ J)(Y)) -
(JY)((\theta\circ J)(X))- \theta([JX,JY])=0, \label{dthetaj}
\end{equation} for all $X,Y$ in $\vfield{TM\slaz}$. For the last equality
in formula \eqref{dthetaj} we used that $\theta\circ J=0$ and
$[JX,JY]=J[X,JY]+J[JX,Y]$, which is true since $N_J=0$. Therefore,
the $2$-form $d\theta$ vanishes on any pair of vertical vectors.
Using the fact that $\nabla J=0$, it follows that the $2$-form
$\nabla d\theta$ also vanishes on any pair of vertical vectors.

For a semi-basic $1$-form $\theta$ we have that $\Phi\circ
J=v\circ \mathcal{L}_S\circ h\circ J=0$ and $J\circ \Phi =J\circ
v\circ \mathcal{L}_S\circ h=0$, since $h\circ J=0$ and $J\circ
v=0$. Therefore,
\begin{equation}
d_{\Phi}\theta(X,JY)= i_{\Phi}d\theta (X,JY)= d\theta(\Phi X,
JY)=0. \label{dphitheta} \end{equation} Last equality in formula
\eqref{dphitheta} is due to the fact that $\Phi X$ and $JY$ are
vertical vector fields.

We evaluate both sides of formula \eqref{nabladt1} on a pair of
vectors of the form $JX, JY$, for arbitrary $X, Y$ in
$\vfield{TM\slaz}$. Using formulae \eqref{dthetaj} and
\eqref{dphitheta} we obtain
\begin{eqnarray}
\label{ifj1} \mathcal{L}_Sd\theta(JX,JY) & = &
i_{\mathbb{F}+J}d\theta(JX, JY)= d\theta\left(hX, JY\right) +
d\theta\left(JX, hY\right) \\ & = & d_J\theta\left(hX,
hY\right)=d_J\theta(X,Y). \nonumber
\end{eqnarray}

We proceed now to prove that $\mathcal{L}_S\theta$ is closed if
and only if conditions \eqref{Helmholtz} are true. From formula
\eqref{nabladt1} it follows that $\mathcal{L}_S \theta$ is closed
if and only if
\begin{eqnarray} \nabla d\theta +
i_{\mathbb{F}+J}d\theta - d_{\Phi}\theta=0. \label{nabladt2}
\end{eqnarray}

We assume first that $\mathcal{L}_S\theta$ is closed and prove
that the four conditions in \eqref{Helmholtz} hold. From formula
\eqref{ifj1} it follows that $d_J\theta=0$. Therefore $\nabla
d_J\theta=\nabla i_Jd\theta=0$. Using the commutation rule $\nabla
i_J=i_J \nabla$, we obtain that $i_J\nabla d\theta=0$ and hence
\begin{eqnarray} (\nabla d\theta)(JX,Y)+(\nabla d\theta)(X,JY)=0,
\forall X, Y \in \vfield{TM\slaz}. \label{eq3}
\end{eqnarray}

Let us evaluate the $2$-form $i_{\mathbb{F}+J}d\theta$ on a pair
of vectors $X, JY$, for $X, Y$ in $\vfield{TM\slaz}$. According to
formula \eqref{ifj1}, this $2$-form vanishes on the pair of
vertical vectors $vX, JY$ and hence we have
\begin{eqnarray}
i_{\mathbb{F}+J}d\theta(X, JY)= i_{\mathbb{F}+J}d\theta(hX, JY) =
d\theta(hX, hY) =d_h\theta(X,Y). \label{ifj2} \end{eqnarray}

Therefore, if we evaluate the left hand side of formula
\eqref{nabladt2} on a pair of vectors $X, JY$, for $X, Y$ in
$\vfield{TM\slaz}$ and use formula \eqref{ifj2} we obtain
\begin{eqnarray}
(\nabla d\theta)(X,JY)+d_h\theta(X,Y)=0. \label{eq5}
\end{eqnarray}
Similarly, if we evaluate the left hand side of formula
\eqref{nabladt2} on a pair of vectors $JX, Y$, for $X, Y$ in
$\vfield{TM\slaz}$ and use formula \eqref{ifj2} we obtain
\begin{eqnarray}
\mathcal{L}_S\theta(X,JY)=(\nabla d\theta)(JX,Y)+d_h\theta(X,Y).
\label{eq6}
\end{eqnarray}

Now, using formulae \eqref{eq5}, \eqref{eq6} and \eqref{eq3} it
follows that $d_h\theta=0$ and $\nabla d\theta=0$. Finally, from
formula \eqref{nabladt2} it follows that last Helmholtz condition
$d_{\Phi}\theta=0$ is also satisfied.

For the other direction, let us assume that conditions in
\eqref{Helmholtz} hold and let us prove that $\mathcal{L}_S\theta$
is closed. In view of formula \eqref{nabladt1}, we only need to
prove that $i_{\mathbb{F}+J}d\theta=0$. Since $(\mathbb{F}+J)\circ
h=0$ it follows that $i_{\mathbb{F}+J}d\theta$ vanishes on any
pair of horizontal vectors. It remains to show that
$i_{\mathbb{F}+J}d\theta(X,JY)=0$, for two arbitrary vector fields
$X$ and $Y$ on $TM\slaz$. For vector field $X$ there exists a
vector field $Z$ on $TM\slaz$ such that $vX=JZ$. Therefore,
\begin{eqnarray*} i_{\mathbb{F}+J}d\theta(X,JY)& = &
i_{\mathbb{F}+J}d\theta(hX,JY) + i_{\mathbb{F}+J}d\theta(JZ,JY)
\\ &=& d_h\theta(X,Y)+ d_J\theta(Z,Y). \end{eqnarray*}
Conditions $d_J\theta=0$, $d_h\theta=0$, and the above
considerations imply that $i_{\mathbb{F}+J}d\theta=0$ and hence
$\mathcal{L}_S\theta$ is closed.
\end{proof}

\subsection{Helmholtz conditions for a multiplier matrix} \label{helm_mm}
We will show how conditions \eqref{Helmholtz}, expressed in terms
of a semi-basic $1$-form, are related with the classic formulation
of Helmholtz conditions expressed in terms of a multiplier matrix.

For a semi-basic $1$-form $\theta=\theta_i dx^i\in
\Lambda^1(TM\slaz)$, let us introduce the following notations
\begin{eqnarray} a_{ij}:=\frac{1}{2}\left(\frac{\delta \theta_i}{\delta
x^j}-\frac{\delta \theta_j}{\delta x^i}\right), \quad
g_{ij}:=\frac{1}{2}\frac{\partial \theta_i}{\partial y^j}.
\label{ab}
\end{eqnarray}
With respect to these notations we have
\begin{eqnarray*}
 d\theta &=& a_{ij} dx^j \wedge dx^i + 2g_{ij} \delta y^j \wedge
 dx^i; \\
 d_h\theta &=& a_{ij} dx^j \wedge dx^i; \\
 d_J\theta &=& (g_{ij}-g_{ji}) dx^j \wedge dx^i; \\
 d_{\Phi}\theta &=& (g_{kj}R^k_i-g_{ik}R^k_j) dx^j \wedge dx^i.
\end{eqnarray*}
Moreover if $d_h\theta=0$ it follows that $\nabla d\theta=2(\nabla
g_{ij})\delta y^j \wedge dx^i$, where $\nabla
g_{ij}=S(g_{ij})-N_i^kg_{kj} - N_j^kg_{ik}$. Therefore, conditions
\eqref{Helmholtz} can be expressed in coordinates as follows:
\begin{eqnarray}
a_{ij}=0, \quad g_{ij}=g_{ji}, \quad \nabla g_{ij}=0,\quad
g_{ik}R^k_j=g_{jk}R^k_i. \label{general_Helmholtz} \end{eqnarray}
Last three conditions in \eqref{general_Helmholtz} together with
$$\frac{\partial g_{ij}}{\partial y^k}= \frac{\partial g_{ik}}{\partial
y^j}, $$ which is satisfied in view second notation \eqref{ab},
are known as the Helmholtz conditions for the inverse problem of
Lagrangian dynamics, \cite{sarlet82}. A global formulation of
Helmholtz conditions \eqref{general_Helmholtz} in terms of the
(0,2)-type symmetric tensor $g=g_{ij}dx^i\otimes dx^j$ along the
tangent bundle projection has been obtained by Martinez et al. in
\cite{martinez93a}.

\subsection{Homogeneous case}

In this section we prove Theorem \ref{exactlst}, which is a
refinement of Theorem \ref{thm_Helm} in the case that the $1$-form
$\theta$ is homogeneous. In this case, $\mathcal{L}_S\theta$ is
closed if and only if $\mathcal{L}_S\theta$ is exact. Also,
depending of the degree of homogeneity, one can drop either one or
two conditions from Helmholtz conditions \eqref{Helmholtz}. See
condition iv) in Theorem \ref{exactlst} below. The fact that for a
spray $S$, one of the Helmholtz condition is a consequence of the
other ones has been proved recently, in a different way, by Prince
\cite{prince08}.

In Proposition \ref{lem_hil} we show that a semi-basic $1$-form,
$(k-1)$-homogeneous with $k\neq 0$, is $d_J$-closed if and only if
it is $d_J$-exact. This result has been obtained in a more general
context by Klein \cite{klein62}, Klein and Voutier \cite{klein68}
and used recently by Vattam\'ani \cite{vattamany04} and Szilasi
and Vattam\'ani \cite{szilasi02} in the Finslerian context.

\begin{prop}\label{lem_hil}
Let $k$ be an integer.
\begin{enumerate}
\item[i)] If $L$ is a $k$-homogeneous function $L\in
  C^\infty(TM\slaz)$, then \emph{Poincar\'e-Cartan} $1$-form $d_JL \in \Lambda^1(TM\slaz)$ is
  semi-basic, $d_J$-closed, and $(k-1)$-homogeneous.
\item[ii)] If a semi-basic $1$-form $\theta\in \Lambda^1(TM\setminus\{0\})$ is
$(k-1)$-homogeneous with $k\neq 0$,  and $d_J$-closed, then
$\theta$ is $d_J$-exact. Moreover, if $S$ is a spray on $M$, then
\begin{eqnarray}
  L=\frac{1}{k}i_S\theta.
\label{listheta}
\end{eqnarray} is the unique $k$-homogeneous function $L\in
  C^\infty(TM\slaz)$ such that $\theta=d_JL$ (we say that $L$ is the \emph{potential function} for the semi-basic $1$-form $\theta$).
  \end{enumerate}
\end{prop}

Let us note that $M$ has at least one spray since we assume that
$M$ is paracompact. Also, by uniqueness in ii), function $L$ in
equation \eqref{listheta} does not depend on $S$.

\begin{proof}
  i) Since the tangent structure $J$ is $0$-homogeneous, which means that $[\setC, J] = -J$, and using
  formula \eqref{dab} we obtain
\begin{eqnarray*}
\mathcal{L}_{\mathbb{C}} d_J L - d_J \mathcal{L}_{\mathbb{C}}L
=-d_JL.
\end{eqnarray*}
  Therefore, $d_JL$ is $(k-1)$-homogeneous since
  $\mathcal{L}_{\mathbb{C}}f=kf$. Also, $d_JL$ is $d_J$-closed by
  equation \eqref{dj2}, and semi-basic since $i_Jd_JL = dL\circ
  J^2=0$.

  ii) Let $S$ be a spray on $M$. We prove that $d_J L=\theta$, when
  function $L$ is defined in equation \eqref{listheta}.  By definition
  we have $JS=\setC$, and by equation \eqref{lsj}, we have $[S,J]= -\Gamma$.
  The generalized Cartan's formula \eqref{ixda} then gives
\begin{eqnarray*}
i_S d_J\theta + d_J i_S\theta = \mathcal{L}_{J(S)}\theta -
i_{[S,J]}\theta = \mathcal{L}_{\mathbb{C}}\theta +
i_{\Gamma}\theta. \label{isdj}
\end{eqnarray*}
Now $d_J\theta=0$, $\mathcal{L}_{\mathbb{C}}\theta=(k-1)\theta$,
and $i_{\Gamma}\theta=\theta\circ \Gamma=\theta\circ h= \theta$,
so $d_J i_S\theta = k \theta$ and $d_J L = \theta$ by equation
\eqref{listheta}.

Let $S$ and $L$ be as in the proof of ii).  To prove that $L$ is
$k$-homogeneous, let us first note that $[\setC, S]=S$, and by
homogeneity $\mathcal{L}_{\mathbb{C}}\theta = (k-1)\theta$.
Commutation rule
\begin{eqnarray*}
  \mathcal{L}_{\mathbb{C}} i_S\theta - i_S \mathcal{L}_{\mathbb{C}} \theta= i_{[\mathbb{C}, S]} \theta
\end{eqnarray*}
then gives $\cL_\setC L = kL$, where $L$ is defined in equation
\eqref{listheta}.

For uniqueness, suppose that $\tilde{L}$ is another
$k$-homogeneous potential function for $\theta$. Then
$\theta=d_J\tilde{L} = d_J L$. If $S^\ast$ is a spray on $M$, then
$d_J L(S^\ast) = \setC(L) = \cL_{\setC} L = k L$. Hence $k L = k
\tilde L$, and $L=\tilde L$.
\end{proof}

\begin{thm}
\label{exactlst} Let $S$ be a spray on $M$, and let $\theta \in
\Lambda^1(TM\setminus\{0\})$ be a semi-basic $1$-form. If $\theta$
is $(k-1)$-homogeneous with $k\neq 0$, then the following
conditions are equivalent:
\begin{itemize}
\item[i)] $\mathcal{L}_S\theta$ is closed;
\item[ii)] $\mathcal{L}_S\theta$ is exact;
\item[iii)] $k\mathcal{L}_S\theta=di_S\theta$;
\item[iv)] $\begin{cases}
                d_h\theta=0,\
                d_J\theta=0,\
                  & \mbox{when} \ k=1, \\
                d_h\theta=0,\
                d_J\theta=0,\
                \nabla d\theta=0,
                  & \mbox{when} \ k\notin \{-1,0,1\}. \\
            \end{cases}$
\end{itemize}
\end{thm}
\begin{proof}
Implications iii) $\Rightarrow$ ii) $\Rightarrow$ i) are clear,
and implication i) $\Rightarrow$ iv) follows by Theorem
\ref{thm_Helm}. To prove implication iv) $\Rightarrow$ iii), let
us assume that one branch in iv) holds. By the generalized
Cartan's formula \eqref{ixda} we have $i_Sd_h\theta + d_hi_S\theta
= \mathcal{L}_S\theta - i_{[S,h]}\theta$. Since $\theta$ is
semi-basic, formula \eqref{jacobiphi} yields
$i_{[S,h]}\theta=\theta \circ \mathbb{F}$, and by assumption
$d_h\theta=0$. Hence
\begin{eqnarray}
   \mathcal{L}_S\theta = d_hi_S\theta + \theta\circ \mathbb{F}.
\label{lst1}
\end{eqnarray}
Since $\theta$ is $d_J$-closed and $(k-1)$-homogeneous,
Proposition \ref{lem_hil} implies that there exists a
$k$-homogeneous function $L\in C^\infty(TM\slaz)$ such that
$kL=i_S\theta$.  Since $J\circ \setF= v$, we have $\theta \circ
\mathbb{F}=d_J f\circ \setF = d_vL$, and using $dL = d_v L + d_h
L$, we obtain
\begin{eqnarray}
   \mathcal{L}_S\theta= kd_hL + d_vL=dL+(k-1)d_hL.
   \label{lst2}
\end{eqnarray}

\noindent \emph{Case $1$:} When $k=1$ equation \eqref{lst2}
implies that $\mathcal{L}_S\theta= dL$ and iii) follows.

\noindent \emph{Case $2$:} We show that if $k\notin \{-1,0,+1\}$
then $d_hL=0$ whence condition iii) follows by equation
\eqref{lst2}.  Using Cartan's formula \eqref{da} we have
$\mathcal{L}_S\theta= i_Sd\theta+di_S\theta= i_Sd\theta+ kdL$.
Combining this with formula \eqref{lst2} gives
$i_Sd\theta=(1-k)d_vL$, whence
$$
  \nabla d_vL = \frac{1}{1-k} \nabla i_Sd\theta=\frac{1}{1-k} i_S\nabla d\theta=0,
$$
where we used $\nabla i_S=i_S\nabla$ and assumption $\nabla
d\theta=0$. Contracting $\nabla d_vL=0$ by $\setC$ similarly gives
$$
  0=i_{\mathbb{C}}\nabla d_vL=\nabla i_{\mathbb{C}} d_vL =k\nabla L.
$$
We have proven that $\nabla L=0$, so $\cL_SL=0$. Equation
\eqref{dab} then gives
\begin{eqnarray}
\label{finalEq}
  d_{[S,J]} L = \cL_S d_J L - d_J \cL_S L = \cL_S d_J L.
\end{eqnarray}
By equations \eqref{hproj}, we have $[S,J] = v-h$, so $d_{[S,J]} L
= d_v L - d_h L$.  Equation \eqref{lst2} gives $\cL_S d_J L =
\cL_S \theta = dL+(k-1)d_hL= d_vL + k d_h L$ since $dL=d_vL +
d_hL$. Now equation \eqref{finalEq} gives $(k+1) d_hL =0$. Thus
$d_hL=0$ and iii) follows.
\end{proof}

\section{The inverse problem of the calculus of variations}

The inverse problem of the calculus of variations for a given
semispray has solutions if and only if there exists a multiplier
matrix that satisfies the Helmholtz conditions
\eqref{general_Helmholtz}, \cite{crampin81, krupkova97, sarlet82}.
Within the Helmholtz conditions, the multiplier matrix is the
Hessian of a (locally defined) Lagrangian for which the given
semispray is an Euler-Lagrange vector field.

In the previous section we did reformulate the Helmholtz
conditions in terms of semi-basic $1$-forms. In this section, we
prove that the inverse problem of the calculus of variation has
solutions if and only if there exists a semi-basic $1$-form that
satisfies the Helmholtz conditions \eqref{Helmholtz}. In this
case, the semi-basic $1$-form is the Poincar\'e-Cartan $1$-form of
a locally defined Lagrangian for which the semispray is an
Euler-Lagrange vector field.

In the homogeneous case, according to Theorem \ref{exactlst}, we
have that if for a spray $S$ there exists a $(k-1)$-homogeneous
semi-basic $1$-form $\theta$, $k\neq 0$, that satisfies the
Helmholtz conditions \eqref{Helmholtz} then its potential function
$L=(1/k)i_S\theta$ is a globally defined Lagrangian for which $S$
is an Euler-Lagrange vector field. We will use this result to
study two inverse problems in Finsler geometry.

\subsection{Lagrangian semisprays}
We show that Helmholtz conditions \eqref{Helmholtz} are necessary
and sufficient conditions for a semispray $S$ to be locally
Lagrangian.

\begin{definition}
\label{def:lagrange} \begin{itemize} \item[i)] A smooth function
$L\in C^{\infty}(TM\slaz)$ is called a \emph{Lagrangian}.
\item[ii)] If for a Lagrangian $L$, its \emph{Poincar\'e-Cartan}
$1$-form $d_JL$ is non-degenerate, then the Lagrangian is called
\emph{regular}. \item[iii)] If there exists a $1$-homogeneous
function $F\in C^{\infty}(TM\slaz)$ such that the Lagrangian
$L=F^2$ is regular, then $F$ is called a \emph{Finsler metric}.
\end{itemize}
\end{definition}
For a regular Lagrangian $L$, the non-degeneracy of the
Poincar\'e-Cartan $1$-form $d_JL$ states that the $n\times n$
symmetric matrix with components
\begin{equation}
g_{ij}(x,y)=\frac{1}{2}\frac{\partial^2 L}{\partial y^i\partial
y^j}(x,y) \label{gij} \end{equation} has rank $n$ on $TM\slaz$,
\cite{morandi90}.

For a Lagrangian $L$, the variational problem leads to the
Euler-Lagrange equations:
\begin{equation}
\frac{d}{dt}\left(\frac{\partial L}{\partial y^i}\right) -
\frac{\partial L}{\partial x^i}=0. \label{el1} \end{equation}

For a semispray $S$, its geodesics, given by the system of second
order differential equations \eqref{eq_geod}, are solutions of the
Euler-Lagrange equations \eqref{el1} if and only if the two sets
of equations are related by formula \eqref{variational_sode}, with
the multiplier matrix given by formula \eqref{gij}. Therefore, if
for a semispray $S$, there exists a Lagrangian $L$ such that
formula \eqref{variational_sode} holds true, then Euler-Lagrange
equations \eqref{el1} are equivalent with \cite{morandi90,
sarlet84}
\begin{equation}
S\left(\frac{\partial L}{\partial y^i}\right) - \frac{\partial
L}{\partial x^i}=0, \label{el3} \end{equation} which can be
further expressed as \begin{equation} \mathcal{L}_Sd_JL=dL.
\label{lsdl}
\end{equation}
For a Lagrangian $L$, a semispray $S$ that satisfies equation
\eqref{lsdl} is called an \emph{Euler-Lagrange vector field}. If
$L$ is regular, $L$ has an unique Euler-Lagrange vector field.

\begin{definition} \label{Slagrangian} A semispray $S$ on
$M$ is called (locally) \emph{Lagrangian} if there exists a
(locally defined) Lagrangian $L$ that satisfies equation
\eqref{lsdl}.\end{definition}

\begin{thm} \label{thm:lagrangian}
Let $S$ be a semispray on $M$. Then, $S$ is a locally Lagrangian
vector field if and only if there exists a semi-basic $1$-form
$\theta \in \Lambda^1(TM\slaz)$ such that the Helmholtz conditions
\eqref{Helmholtz} are satisfied.
\end{thm}
\begin{proof}
We assume that the semispray $S$ is derived from a locally defined
Lagrangian $L$. Consider $\theta=d_JL$, the Poincar\'e-Cartan
$1$-form of $L$. From formula \eqref{lsdl} it follows that
$\mathcal{L}_S\theta$ is closed and using Theorem \ref{thm_Helm}
it follows that the semi-basic $1$-form $\theta$ satisfies
Helmholtz conditions \eqref{Helmholtz}.

For the converse, consider a semi-basic $1$-form $\theta\in
\Lambda^1(TM\slaz)$ such that Helmholtz conditions
\eqref{Helmholtz} are satisfied. Using Theorem \ref{thm_Helm} it
follows that the $1$-form $\mathcal{L}_S\theta$ is closed.
Therefore, there exists a locally defined function $L$ on
$TM\slaz$ such that
\begin{equation} \mathcal{L}_S\theta=dL. \label{eq:lstheta}
\end{equation} If we apply $i_J$ to both sides of formula
\eqref{eq:lstheta} we obtain
\begin{equation} i_J\mathcal{L}_S\theta=d_JL. \label{ijls}
\end{equation} From formulae \eqref{lxia} and \eqref{lsj} we
obtain the following commutation formula
\begin{equation}
i_J\mathcal{L}_S - \mathcal{L}_S i_J = - i_{[S,J]}=i_{h-v}.
\label{com_ijls} \end{equation} Now, we substitute the derivation
$i_J\mathcal{L}_S$ from formula \eqref{com_ijls} into formula
\eqref{ijls}, we use that $\theta$ is semi-basic, which implies
that  $i_J\theta=0$ and $i_{h-v}\theta=\theta$ and obtain
\begin{equation} d_JL=\theta. \label{djl} \end{equation} In
view of equations \eqref{eq:lstheta} and \eqref{djl} we obtain
that equation \eqref{lsdl} is satisfied and hence the semispray
$S$ is a locally Lagrangian vector field.
\end{proof}
The regularity of a Lagrangian is characterized by the
non-degeneracy of its Poincar\'e-Cartan $1$-form. Therefore, a
semispray $S$ is induced by a (locally defined) regular Lagrangian
if and only if there exists a non-degenerate semi-basic $1$-form
$\theta \in \Lambda^1(TM\slaz)$ that satisfies the Helmholtz
conditions \eqref{Helmholtz}.

Theorem \ref{thm:lagrangian} was inspired by a Theorem of Crampin
\cite{crampin81}, where locally Lagrangian semisprays are
characterized in terms of $2$-forms. A version of this result, in
the homogeneous case, is due to Klein \cite{klein62}.

Sarlet et al. \cite{sarlet84} associate to a semispray $S$ a
particular subset
$\Lambda^1_S(TM\slaz)=\{\omega\in\Lambda^1(TM\slaz),
\mathcal{L}_Si_J\omega=\omega\}$ of $1$-forms on $TM\slaz$.
(Locally) Lagrangian semisprays are then characterized by the
property that $\Lambda^1_S(TM\slaz)$ contains an element $\omega$
that is (closed) exact and $i_J\omega$ is non-degenerate. The
relation between this result and Theorem \ref{thm:lagrangian} is
as follows. Let $\theta$ be a non-degenerate, semi-basic $1$-form
such that $\mathcal{L}_S\theta$ is closed. Consider the closed
$1$-form $\omega=\mathcal{L}_S\theta$. From formula \eqref{ijls}
it follows that $i_J\omega=\theta$ is non-degenerate, and from
equation \eqref{lsdl} it follows that
$\mathcal{L}_Si_J\omega=\omega$, which means that $\omega \in
\Lambda^1_S(TM\slaz)$.

For a Lagrangian semispray $S$, two of the Helmholtz conditions
\eqref{general_Helmholtz}: $a_{ij}=0$ and $\nabla g_{ij}=0$ where
used in \cite{bucataru07} to characterize the canonical nonlinear
connection of a Lagrange space.

\subsection{Further discussions of Helmholtz conditions} For a
semispray $S$, consider a semi-basic $1$-form $\theta$ on
$TM\slaz$ that satisfies the Helmholtz conditions
\eqref{Helmholtz}. Three of these conditions can be expressed as
follows
\begin{eqnarray}
i_{\Gamma}d\theta=0, \quad i_{J}d\theta=0, \quad
i_{\Phi}d\theta=0. \label{Halgebric} \end{eqnarray} First two
conditions \eqref{Halgebric} fixes the number of unknown
components of $d\theta=2g_{ij}\delta y^j \wedge dx^i$ to
$n(n+1)/2$. Third condition \eqref{Halgebric} imposes algebraic
restrictions on $d\theta$.

Grifone and Muzsnay associate to a semispray $S$ the graded Lie
algebra $\mathcal{A}_S$ of vector valued forms $A$ such that
$i_Ad\theta=0$. Using Theorem \ref{prop:nabla} it follows that if
$A\in \mathcal{A}_S$ then $\nabla A \in \mathcal{A}_S$. Therefore,
iterated covariant derivatives $\nabla^k\Phi$ of the Jacobi
endomorphism impose further algebraic restrictions on $d\theta$
\begin{eqnarray}
i_{\nabla^k \Phi}d\theta=0. \label{inkp} \end{eqnarray} The
sequence of $(1,1)$-type tensor fields $\Phi^{(k)}:=\nabla^k\Phi$
where considered previously by Sarlet \cite{sarlet82}, Crampin
\cite{crampin83} and Grifone and Muzsnay \cite{grifone00}.

From formula \eqref{dab} it follows that if $A, B\in
\mathcal{A}_S$ then $[A,B]\in \mathcal{A}_S$. Therefore, Helmholtz
conditions $d_J\theta=0$ and $d_{\Phi}\theta=0$ and formula
\eqref{jphir} imply that $d_R\theta=0$, which gives a new
algebraic restriction on $d\theta$
\begin{eqnarray}
i_{R}d\theta=0. \label{irdt} \end{eqnarray} Hence, the graded Lie
algebra $\mathcal{A}_S$ of algebraic restrictions on $d\theta$
contains also the sequence of iterated covariant derivatives
$\nabla^kR$ of the curvature tensor $R$.

The graded Lie algebra $\mathcal{A}_S$ is used in general to
formulate non-existence results for a semispray $S$ to be
Lagrangian, \cite{grifone00,sarlet82}. It follows that if there
exists $p\in M$ such that
$\operatorname{rank}\{\mathcal{A}_S(p)\}>n(n+1)/2$ then $S$ is not
Lagrangian.

We note that for the homogeneous case the fact that some of the
Helmholtz conditions can be derived from the other ones, in a
non-linear way, does not change the rank of $\mathcal{A}_S$ and
hence it does not change the rank of algebraic restrictions one
have to impose on $d\theta$.

\subsection{Lagrangian sprays} \label{lagrangianspray}

We show that in the homogeneous case, a spray $S$ is Lagrangian if
and only if only two or three of the Helmholtz conditions are
satisfied, depending on the degree of homogeneity. In particular
we discuss Helmholtz conditions for two important inverse
problems: projective metrizability and Finsler metrizability.

\begin{thm}
\label{thm_homogeneous} Let $S$ be a spray on $M$. Then $S$ is a
Lagrangian vector field, induced by a $k$-homogeneous Lagrangian,
if and only if there exists a $(k-1)$-homogeneous, semi-basic
$1$-form $\theta \in \Lambda^1(TM\setminus\{0\})$ such that
\begin{equation}\begin{cases}
                d_h\theta=0,\
                d_J\theta=0,\
                  & \mbox{when} \ k=1, \\
                d_h\theta=0,\
                d_J\theta=0,\
                \nabla d\theta=0,
                  & \mbox{when} \ k\notin \{-1,0,1\}. \\
            \end{cases} \label{helm_hom}
\end{equation}
\end{thm}
\begin{proof}
Suppose that the spray $S$ is an Euler-Lagrange vector field for a
$k$-homogenous lagrangian $L$. It follows that the
Poincar\'e-Cartan $1$-form $\theta=d_JL$ is a $(k-1)$-homogeneous,
semi-basic $1$-form. Since equation \eqref{lsdl} holds true it
follows using Theorem \ref{exactlst} that $\theta$ satisfies
conditions \eqref{helm_hom}.

Conversely, suppose that there exists a $(k-1)$-homogeneous,
semi-basic $1$-form $\theta \in \Lambda^1(TM\setminus\{0\})$
satisfies conditions \eqref{helm_hom}. From Proposition
\ref{lem_hil} it follows that $L=(1/k)i_S\theta$ is a
$k$-homogeneous Lagrangian. Using Theorem \ref{exactlst} it
follows that conditions \eqref{helm_hom} imply that $L$ satisfies
equation \eqref{lsdl} and hence $S$ is a Lagrangian vector field.
\end{proof}
Next, we will discuss in more details the two branches of
conditions \eqref{helm_hom} and show that they correspond to two
inverse problems studied in Finsler geometry: Finsler
metrizability and projective metrizability.

\begin{definition} A spray $S$ is \emph{projectively metrizable}
if there exists a $1$-homogeneous Lagrangian $F$ such that
equation \eqref{lsdl} is satisfied. \end{definition} Note that in
this definition and hence this work we do not necessarily assume
that $F$ is a Finsler metric, which in addition would require that
the Hessian of $F$ with respect to the fibre derivatives has rank
$(n-1)$. For a discussion on the regularity of the Lagrangian
$L=F^2$ and the hessian of $F$ we refer to the book of Matsumoto
\cite{matsumoto86} as well as to the recent work of Crampin
\cite{crampin08} and Szilasi \cite{szilasi07}. If a spray $S$ is
projectively metrizable, its geodesics, up to an orientation
preserving reparameterization, are solutions of the Euler-Lagrange
equations of a $1$-homogeneous lagrangian $L$. Indeed if a $F$ is
a $1$-homogenous solution of \eqref{lsdl} then the Euler-Lagrange
equations \eqref{el1} for $F$ can be written as
\begin{eqnarray}
h_{ij}\left(x,\frac{dx}{dt}\right)\left(\frac{d^2x^i}{dt^2}+
2G^i\left(x,\frac{dx}{dt}\right)\right)=0.\label{hijg}\end{eqnarray}
In the above equations \eqref{hijg} $h_{ij}$ are the components of
the Hessian of $F$ with respect to the fiber coordinates. Since
$h_{ij}$ are $(-1)$-homogeneous it follows that
$h_{ij}\frac{dx^j}{dt}=0$ and hence the system of equations
\eqref{hijg} is invariant under an orientation preserving
reparameterization.

The problem of projective metrizability of a spray $S$ is related
to Hilbert's fourth problem. For a flat spray this problem was
first studied by Hamel \cite{hamel03} and it is known as the
Finslerian version of Hilbert's fourth problem \cite{alvarez05,
crampin08, szilasi07}. For a general spray, Rapcs\'ak
\cite{rapcsak62} was first to provide criteria, in local
coordinates, for the projective metrizability of a spray. Global
formulations for projective metrizability criteria were obtained
by Klein \cite{klein62}, Klein and Voutier \cite{klein68} and
Szilasi \cite{szilasi07}. An extensive discussion of the
projective metrizability of a spray appears in Vattam\'any's Ph.D
thesis \cite[chapter 2]{vattamany04}.

According to Theorem \ref{thm_homogeneous}, we have that a spray
$S$ is projectively metrizable if and only if there exists a
$0$-homogeneous, semi-basic $1$-form $\theta\in
\Lambda^1(TM\slaz)$ such that $d_h\theta=0$ and $d_J\theta=0$.
According to Proposition \ref{lem_hil} the condition $d_J\theta=0$
implies that $F=i_S\theta$ is the only $1$-homogeneous Lagrangian
that satisfies $\theta=d_JF$. Moreover, from Theorem
\ref{exactlst} it follows that $F$ satisfies the condition
$\mathcal{L}_Sd_JF=dF$, which is equivalent to $i_Sdd_JF=0$. Last
condition represents condition Rap 1 in Theorem 8.1 by Szilasi
\cite{szilasi07}. Also condition $d_h\theta=0$ represents
condition Rap 4 in the same cited work.

For the particular case of a flat spray we obtain that the induced
nonlinear connection is integrable and hence $[h,h]=0$. It follows
that $d_h^2=0$ and therefore any $d_h$-closed semi-basic $1$-form
is locally $d_h$-exact, \cite{vaisman73}. Since $d_h\theta=0$ it
follows that there exists a $0$-homogeneous function $f\in
C^{\infty}(TM\slaz)$ such that $\theta=d_hf$. From the above
discussion we have that the $1$-homogeneous function $F=i_S\theta$
projectively metricizes the spray $S$ if and only if
$\theta=d_hf$. Therefore
$$ F=i_S\theta=i_Sd_hf=S(f)$$ projectively metricizes the spray
$S$ if and only if $d_hd_Jf=0$. In local coordinates, we have that
last condition is equivalent to \begin{equation} \frac{\partial^2
f}{\partial y^i\partial x^j}= \frac{\partial^2 f}{\partial
y^j\partial x^i}.\label{f2xy}\end{equation} This is a
reformulation of Proposition 2 by Crampin \cite{crampin08} or
Proposition 8.1 by Szilasi \cite{szilasi07}, which state that
$F=S(f)$ is a $1$-homogeneous function that projectively
metricizes the spray $S$ if and only if there exists a
$0$-homogeneous function $f$ on $TM\slaz$ that satisfies condition
\eqref{f2xy}. Both Crampin and Sarlet ask more conditions for the
symmetric bilinear form with components \eqref{f2xy} to obtain
that $F=S(f)$ is a Finsler function.

\begin{definition} A spray $S$ is \emph{Finsler metrizable}
if there exists a $2$-homogeneous Lagrangian $L$ such that
equation \eqref{lsdl} is satisfied. \end{definition} Note that in
this definition and hence this work we do not necessarily require
the regularity of the Lagrangian. If a spray $S$ is Finsler
metrizable, its geodesics are also solutions of the Euler-Lagrange
equations of a $2$-homogeneous lagrangian $L$. The Finsler
metrizability problem, viewed as the inverse problem of the
calculus of variation restricted to the class of $2$-homogeneous
Lagrangians has been studied recently by Crampin \cite{crampin07},
Krupka and Sattarov \cite{krupka85}, Muzsnay \cite{muzsnay06},
Prince \cite{prince08}, Szilasi and Vattam\'ani \cite{szilasi02}.

According to Theorem \ref{thm_homogeneous}, we have that a spray
$S$ is Finsler metrizable if and only if there exists a
$1$-homogeneous, semi-basic $1$-form $\theta\in
\Lambda^1(TM\slaz)$ such that $d_h\theta=0$, $d_J\theta=0$ and
$\nabla d\theta=0$. According to Proposition \ref{lem_hil} the
condition $d_J\theta=0$ implies that $2L=i_S\theta$ is the only
$2$-homogeneous Lagrangian that satisfies $\theta=d_JL$. Moreover,
from Theorem \ref{exactlst} it follows that $L$ satisfies the
condition $\mathcal{L}_Sd_JL=dL$, which is equivalent to
$i_Sdd_JL=-dL$. Last condition is equivalent to $d_hL=0$ that has
been used by Muzsnay \cite{muzsnay06} to obtain necessary and
sufficient conditions for Finsler metrizability in term of an
associated holonomy algebra.

\begin{acknowledgement*}
The authors wish to acknowledge fruitful comments from W. Sarlet
and express their thanks to J. Szilasi for the reference
\cite{vattamany04}.

I.B has been supported by grant ID 398 from the Romanian Ministry
of Education. M.D. has been supported by Academy of Finland Center
of Excellence Programme 213476, the Institute of Mathematics at
the Helsinki University of Technology, and Tekes project MASIT03
-- Inverse Problems and Reliability of Models.
\end{acknowledgement*}

\end{document}